\documentclass{article}
\usepackage{graphicx}
\usepackage{amsmath}
\usepackage[symbol]{footmisc}
\usepackage[french]{babel}

\newtheorem{theorem}{Theorem}[section]

\newtheorem{remark}[theorem]{Remark}

\newenvironment{proof}[1][Proof]{\textbf{#1.} }{\ \rule{0.5em}{0.5em}}
\input{tcilatex}
\renewcommand{\theequation}{\thesection.\arabic{equation}}

\begin{document}

\begin{center}
{\huge A recurrence for certain Tutte polynomials}

\bigskip

{\large Vincent Brugidou} 
\footnotetext{\textit{E-mail address:}vincent.brugidou@univ-lille.fr}
\end{center}

\bigskip

\begin{center}
\textit{Universit\'{e} de Lille, 59655 Villeneuve d'Ascq cedex, France}
\end{center}

\bigskip

\textbf{Abstract }We combinatorially prove a new recurrence between the
Tutte polynomials of graphs obtained by contraction of the complete graphs $%
K_{n}$. This generalizes, to two variables, a relation previously obtained
by the author between the inversion enumerator polynomials in certain
colored tree sequences.

\bigskip

\textit{keywords: }Tutte polynomial, complete graphs, contracted graphs,
tree inversions.

.

\section{\protect\bigskip \textbf{Introduction}}

\bigskip

In [2, Eq. (6.5)], we have algebraically proven the following recurrence
relation valid for integers $n$ and $r$, such that $n>r\geq 1$:

\begin{equation}
J_{n}^{\left( r\right) }(q)=\sum\limits_{s=1}^{n-r}\binom{n-r}{s}\left[ r%
\right] _{q}^{s}\;q^{\binom{s}{2}}J_{n-r}^{\left( s\right) }(q).  \tag{1.1}
\end{equation}
In this equation, $\left[ r\right] _{q}$ is the $q$-analogue of $r$ and the
polynomials $J_{n}^{\left( r\right) }$ are special cases of the inversion
enumerator polynomials in colored tree sequences introduced by Stanley [7]
and Yan [8]. These latter polynomials are themselves generalizations of the
polynomial $J_{n}=J_{n}^{\left( 1\right) }$, which is the enumerator of
inversions in trees on $n$ vertices, first introduced in [6]. All these
polynomials also have interpretations in terms of parking functions, for
which we refer to the review article [9].

In [4], it was noted that: 
\begin{equation}
J_{n}^{\left( r\right) }(q)=T_{n}^{(r)}\left( 1,q\right) \text{,}  \tag{1.2}
\end{equation}

where $T_{n}^{\left( r\right) }\left( x,y\right) $ is the Tutte polynomial
of the complete graph $K_{n}$, contracted over all edges between $r$
vertices of $K_{n}$. We will sometimes use the simplified notation $K_{n/r}$
for this contracted graph (see [3, Remark 4.5] for the explanation of this
notation).

The main result of this article is the following theorem, which generalizes $%
\left( 1.1\right) $ to two variables.

\begin{theorem}
The Tutte polynomial of graphs $K_{n/r}$ verifies for $n>r\geq 1$ the
recurrence relation: 
\begin{equation}
T_{n}^{\left( r\right) }(x,y)=\sum\limits_{s=1}^{n-r}\binom{n-r}{s}\left[ r%
\right] _{y}^{s}\;y^{\binom{s}{2}}T_{n-r}^{\left( s\right) }(x,y)+\left(
x-1\right) T_{n-r}^{\left( 1\right) }\left( x,y\right) \text{,}  \tag{1.3}
\end{equation}
where $\left[ r\right] _{y}=1+y+y^{2}+...+y^{r-1}$.
\end{theorem}

This equation is proven combinatorially, using one of the classical
expressions for the Tutte polynomial. It allows us to recursively compute
all polynomials $T_{n}^{\left( r\right) }$ (see Table 1), and in particular
the Tutte polynomials of complete graphs $K_{n}$ since $T_{K_{n}}=T_{n}^{%
\left( 1\right) }$.

By restricting to connected graphs, we show in Section 4 how a combinatorial
proof of Equation (1.1) can be easily obtained as a special case of the
proof of Theorem 1.1.

\section{Notation Conventions}

\subsection{$\protect\bigskip $On Sets}

If $A$ and $B$ are finite sets, $\left| A\right| $ denotes the cardinality
of $A$ and $A-B=\left\{ x\in A,x\notin B\right\} $.

\subsection{$\protect\bigskip $On Graphs}

The terminology and notations concerning graphs are, in general, those of
[1]. We recall here some of these notations that are important for our
subject.

Let $G=\left( V\left( G\right) ,E\left( G\right) \right) $ be a finite graph
whose vertex set is $V\left( G\right) $ and whose edge set is $E\left(
G\right) $. If $F\subseteq $ $E\left( G\right) $, then $G\backslash F$ is
the graph with the same vertices as $G$ and whose edge set is $E\left(
G\right) -F$. Still with $F\subseteq $ $E\left( G\right) $, $G/F$ is the
graph obtained by contraction of the edges belonging to $F$. If $R\subseteq
V\left( G\right) $, then $G-R$ is the graph obtained by deleting from $G$
the vertices belonging to $R$, along with all edges incident to those
vertices. If $R$ contains no vertices adjacent in $G$, then $G/R$ is the
graph obtained by replacing all the vertices in $R$ with a single vertex,
which is incident to all the edges incident to the elements of $R$.

We also adopt the author's own notational conventions, as follows. If $V$ is
a finite set, we denote by $K_{V}$ the complete graph whose vertex set is $V$%
. Clearly, $K_{V}\approx K_{\left| V\right| }$. Let $V\mathbf{=}\left\{
1,2,...,n\right\} $, and $R=\left\{ u_{1},u_{2},...,u_{r}\right\} \subseteq
V $ with $\left| R\right| =r$. Define $\overline{R}=V-R=\left\{
w_{1},w_{2},...,w_{n-r}\right\} $. In $K_{V}/E\left( K_{R}\right) \approx
K_{n/r}$, we denote by $0_{R}$ the vertex that replaces, in the contraction,
the $r$ elements of $R$, and by $(u_{i},w_{k})$ the edge incident to $0_{R}$%
, coming under the contraction from the edge joining $u_{i}$ to\ $w_{k}$ in $%
K_{V}$. The edges joigning two vertices of $\overline{R}$ remain unchanged
under contraction.

Let $S=\left\{ v_{1},v_{2},...v_{s}\right\} \subseteq \overline{R}$ with $%
\left| S\right| =s$, and contract $K_{V}/E\left( K_{R}\right) $ along the
edges of $K_{S}$. We obtain the set: 
\begin{equation*}
(K_{V}/E\left( K_{r}\right) )/E\left( K_{S}\right) =K_{V}/(E(K_{R})\cup
E\left( K_{S}\right) ).
\end{equation*}
The previous edge-notation procedure can be applied to this graph: $0_{S}$
will denote the vertex replacing the vertices belonging to $S$, an edge
joigning $0_{S}$ to $w\in (\overline{R}-S)$ will be denoted $(v_{j},w)$ and
an edge joigning $0_{R}$ to $0_{S}$ will be denoted $(u_{i},v_{j})$. The
others edges remain unchanged with respect to those of $K_{V}/E\left(
K_{R}\right) $.

Figure 1 illustrates the notations thus obtained for the edges adjaccent to $%
0_{R}$ or $0_{S}$, in the graphs $K_{5}/E\left( K_{R}\right) $ and $%
K_{5}/(E(K_{R})\cup E(K_{S}))$, for the case where $n=5$, $R=\left\{
1,2\right\} $ and $S=\left\{ 3,4\right\} $.

\begin{figure}[!htp]
\begin{center}
\includegraphics[width=12cm]{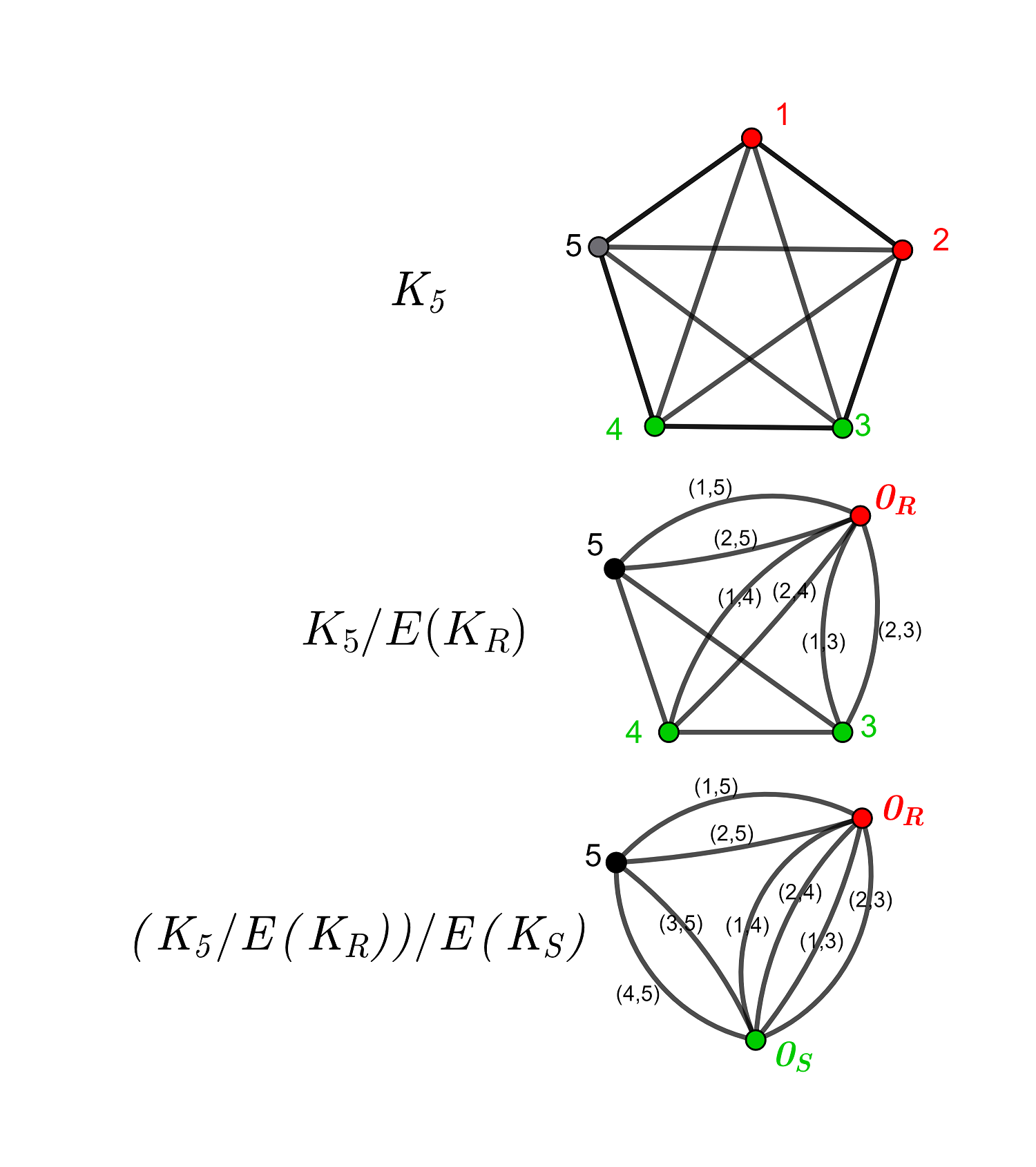}
\end{center}
\caption{Notations obtained for the edges in the graphs $K_5/E(K_R)$ and $%
(K_5/E(K_R))/E(K_S)$, for $n=5$, $R=\{1,2\}$ and $S=\{3,4\}$. The vertices
belonging to $R$ and $S$ are colored red and green, respectively.}
\end{figure}

\section{Proof of Theorem 1}

Let $G=(V\left( G\right) ,E\left( G\right) )$ be a finite graph. We denote
by $c\left( G\right) $ the number of connected components of $G$, by $%
e\left( G\right) =\left| E\left( G\right) \right| $ the number of edges of $G
$, and by $v\left( G\right) =\left| V\left( G\right) \right| $ the number of
vertices of $G$.

\begin{proof}
We start from the following definition of the Tutte polynomial of $K_{n/r}$
(see, for example, Wikipedia, \textit{Tutte polynomial}): 
\begin{equation}
T_{n}^{\left( r\right) }\left( x,y\right) =\sum\limits_{H}\left( x-1\right)
^{c\left( H\right) -1}\left( y-1\right) ^{e\left( H\right) +c\left( H\right)
-v\left( H\right) }  \tag{3.1}
\end{equation}
where the sum runs over the set of spanning subgraphs $H$ of $K_{n/r}\approx
K_{V}/E\left( K_{R}\right) $, with $V=\left\{ 1,2,...,n\right\} $ and $%
R=\left\{ u_{1},u_{2},...,u_{r}\right\} $. This set will be denote by $%
\mathcal{H}(K_{n/r})$. Let $S$ be any subset of $\overline{R}$, and let $%
\mathcal{H}(K_{n/r},S)$ be the set of spanning subgraphs of $K_{n/r}$ such
that the set of vertices adjacent to $0_{R}$ is exactly $S$. We clearly
have: 
\begin{equation}
\mathcal{H}\left( K_{n/r}\right) =\biguplus\limits_{S\subseteq \overline{R}}%
\mathcal{H}\left( K_{n/r},S\right) =\biguplus_{s=0}^{n-r}\biguplus_{\left|
S\right| =s,S\subseteq \overline{R}}\mathcal{H}\left( K_{n/r},S\right) \text{%
.}  \tag{3.2}
\end{equation}
And therefore 
\begin{equation*}
T_{n}^{\left( r\right) }\left( x,y\right) =\sum\limits_{S\subseteq \overline{%
R}}T_{S}\left( x,y\right) 
\end{equation*}
by defining 
\begin{equation}
T_{S}\left( x,y\right) =\sum\limits_{H\in \mathcal{H}\left( K_{n/r},S\right)
}\left( x-1\right) ^{c\left( H\right) -1}\left( y-1\right) ^{e\left(
H\right) +c\left( H\right) -v\left( H\right) }\text{.}  \tag{3.3}
\end{equation}
We will calculate $T_{S}\left( x,y\right) $ in three-steps depending on the
value of $s=\left| S\right| $.\hspace{0.14in}

\hspace{0.14in}

1) $s=0\Leftrightarrow S=\emptyset $

Consider a graph $H\in \mathcal{H}\left( K_{n/r},\emptyset \right) $ and
define $H^{\prime }=\theta \left( H\right) =H-0_{R}$. This defines a
bijection $\theta $ from $\mathcal{H}\left( K_{n/r},\emptyset \right) $ to
the set of spanning subgraphs of $K_{V-R}\approx K_{n-r}$, whith the inverse
bijection given by $\theta ^{-1}\left( H^{\prime }\right) =H^{\prime }+0_{R}$%
. Here, $H^{\prime }+0_{R}$ is the graph whose vertex set is $V\left(
H^{\prime }\right) \cup \left\{ 0_{R}\right\} $ and whose edges are the same
as those of $H^{\prime }$. Between $H$ and $H^{\prime }$, there are the
relations:

\begin{equation*}
c\left( H\right) =c\left( H^{\prime }\right) +1,\;e\left( H\right) =e\left(
H^{\prime }\right) ,\;v\left( H\right) =v\left( H^{\prime }\right) +1\text{.}
\end{equation*}
Consequently, we have: 
\begin{equation*}
T_{\emptyset }\left( x,y\right) =\left( x-1\right) \sum\limits_{H^{\prime
}}\left( x-1\right) ^{c\left( H^{\prime }\right) -1}\left( y-1\right)
^{e\left( H^{\prime }\right) +c\left( H^{\prime }\right) -v\left( H^{\prime
}\right) }=\left( x-1\right) T_{n-r}\left( x,y\right) \text{,}
\end{equation*}
where, in the sum above, $H^{\prime }$ runs over the spanning subgraphs of $%
K_{n-r}$.

\hspace{0.14in}

2) s=1$\Leftrightarrow S=\left\{ v_{1}\right\} $.

There are $\binom{n-r}{1}$ possible choices for the vertex $v_{1}$ $\in $ $%
\overline{R}$. Fix $v_{1}$, and let $H\in \mathcal{H}\left( K_{n/r},\left\{
v_{1}\right\} \right) $. Consider the map $\theta $\ defined by $H^{\prime
}=\theta \left( H\right) =H-0_{R}$. Then $H^{\prime }$ is a spanning
subgraph of $K_{n-r}$. Conversely, if $H^{\prime }$ is a spanning subgraph
of $K_{n-r}$, its inverse image $\theta ^{-1}\left( H^{\prime }\right) $ is
obtained by adding $l$ edges among the $r$ multiple labeled edges $%
(u_{i},v_{1})$ with $u_{i}\in R$. The integer $l$ is between $1$ and $r$,
because by definition of $\mathcal{H}\left( K_{n/r},\left\{ v_{1}\right\}
\right) $, there is at least one edge connecting $v_{1}$ to $0_{R}$. Once $l$
is fixed, there are $\binom{r}{l}$ possible choices for these $l$ edges. For
each choice, denoting the resulting graph by $H$, we have: 
\begin{equation*}
v\left( H\right) =v\left( H^{\prime }\right) +1,\;e\left( H\right) =e\left(
H^{\prime }\right) +l,\;\;c\left( H\right) =c\left( H^{\prime }\right) \text{%
.}
\end{equation*}
The last equality holds because $v_{1}$ is always connected to $0_{R}$, \
and all other edges remain the same in $H^{\prime }$ and $H$. Therefore,
when summing over all possible choice of the $l$ edges, we have: 
\begin{equation*}
T_{\left\{ v_{1}\right\} }\left( x,y\right) =\sum\limits_{l=1}^{r}\binom{r}{l%
}\left( y-1\right) ^{l-1}\sum\limits_{H^{\prime }}\left( x-1\right)
^{c\left( H^{\prime }\right) -1}\left( y-1\right) ^{e\left( H^{\prime
}\right) +c\left( H^{\prime }\right) -v\left( H^{\prime }\right) }\text{,}
\end{equation*}
where the last sum runs over all spanning subgraphs of $K_{n-r}$. Since 
\begin{equation}
\sum\limits_{l=1}^{r}\binom{r}{l}\left( y-1\right) ^{l-1}=\left( y-1\right)
^{-1}\left( \sum\limits_{l=0}^{r}\binom{r}{l}\left( y-1\right) ^{l}-1\right)
=\dfrac{y^{r}-1}{y-1}=\left[ r\right] _{y}\text{,}  \tag{3.4}
\end{equation}
we obtain the following contribution for all $S$ such that $\left| S\right|
=1$:

\begin{equation*}
\binom{n-r}{1}\left[ r\right] _{y}T_{n-r}(x,y)\text{.}
\end{equation*}

\hspace{0.14in}

3) $2\leq s\leq n-r$.

There are $\binom{n-r}{s}$ possible choices of $S$ $\subseteq $ $\overline{R}
$ with $\left| S\right| =s$. Fix $S$, and let $S=\left\{
v_{1},v_{2},...,v_{s}\right\} $ and $\overline{R}-S=\left\{
w_{1,}w_{2},...,w_{n-r-s}\right\} $. Let $H\in \mathcal{H}\left(
K_{n/r},S\right) $, and consider the map $\phi $ defined by 
\begin{equation*}
\phi \left( H\right) =H^{\prime }=\left( H-0_{R}\right) /E\left(
K_{S}\right) \text{.}
\end{equation*}
$H^{\prime }$ is a spanning subgraph of $K_{V-R}/E\left( K_{S}\right)
\approx K_{n-r/s}$. Conversely, given $H^{\prime }$ as\ a spanning subgraph
of $K_{n-r/s}$, we are going to enumarate all graphs in $\mathcal{H}\left(
K_{n/r},S\right) $ that belong to $\phi ^{-1}\left( H^{\prime }\right) $. To
do this, we will describe $\phi $ as a composition of three maps.

First, we have $\phi =\gamma _{1}\circ \theta _{1}$, where $\theta _{1}$ is
the map defined by $\theta _{1}\left( H\right) =H-0_{R}$ and $\gamma _{1}$
is the map \ that, to a graph belonging to $\theta _{1}\left( \mathcal{H}%
\left( K_{n/r},S\right) \right) $, associates its contracted along the edges
of $K_{S}$. Note that $\theta _{1}$ and $\gamma _{1}$ act respectively on
the vertex subsets $R$ and $S$, which are disjoint. We can therefore commute
their action by writing: 
\begin{equation*}
\phi \left( H\right) =H/E\left( K_{S}\right) -0_{R}=\theta _{2}\circ \gamma
_{2}\left( H\right) \text{,}
\end{equation*}
where $\gamma _{2}\left( H\right) =H^{0}=H/E\left( K_{S}\right) $ and $%
\theta _{2}\left( H^{0}\right) =H^{0}-0_{R}$. Since $H$ $\ $and $H^{0}$ are
spanning graphs of $K_{n}/E\left( K_{r}\right) $ and $K_{n}/(E\left(
K_{R}\right) \cup E\left( K_{S}\right) )$, respectively, we will adopt the
notations introduced in Section 2 to denote their edges.

Returning to the definitions recalled in Section 2, we observe that $\gamma
_{2}=\beta _{2}\circ \alpha _{2}$ with:

- $\alpha _{2}\left( H\right) =H^{\ast }=H\backslash E\left( K_{S}\right) $,
that is, $\alpha _{2}$ removes in $H$ all edges joigning two vertices in $S.$
Thus, in $H^{\ast }$, no two element of $S$\ are adjacent.

- $\beta _{2}\left( H^{\ast }\right) =H^{0}=H^{\ast }/S$. Consequently, all
the elements of $S$\ are merged into a single vertex $0_{S}$, and the edges $%
(u_{i},v_{j})$ and $(v_{j},w_{k})$ of $H^{\ast }$ become edges with the same
label in $H^{0}$, while the other edges remain unchanged. It is clear that $%
\beta _{2}$ is a bijection from $\alpha _{2}\left( \mathcal{H}\left(
K_{n/r},S\right) \right) $ into $\gamma _{2}\left( \mathcal{H}\left(
K_{n/r},S\right) \right) $, whose reciprocal bijection $\beta _{2}^{-1}$ is
obtained by splitting $0_{S}$ into the $s$ elements of $S$. We note that $%
H^{\ast }$ and $H^{0}$ have the same number of connected components, since
the elements of $S$ in $H^{\ast }$ are in the same connected component
containing $0_{R}$. Thus, we have: 
\begin{equation}
v\left( H^{\ast }\right) =v\left( H^{0}\right) +s-1,\;\;e\left( H^{\ast
}\right) =e\left( H^{0}\right) ,\;\;c(H^{\ast })=c\left( H^{0}\right) \text{.%
}  \tag{3.5}
\end{equation}
We therefore finally have the decomposition $\phi \left( H\right) =\theta
_{2}\circ \beta _{2}\circ \alpha _{2}\left( H\right) =H^{\prime }$, hence $%
\phi ^{-1}\left( H^{\prime }\right) =\alpha _{2}^{-1}\left( \beta
_{2}^{-1}\left( \theta _{2}^{-1}\left( H^{\prime }\right) \right) \right) $.
In Figure 2, we have represented the successive images, under these maps, of
a graph $H\in \mathcal{H}\left( K_{n/r},S\right) $ with $n=8$, $R=\left\{
1,2\right\} $ and $S=\left\{ 3,4,6\right\} $.

\begin{figure}[tph]
\begin{center}
\includegraphics[width=12cm]{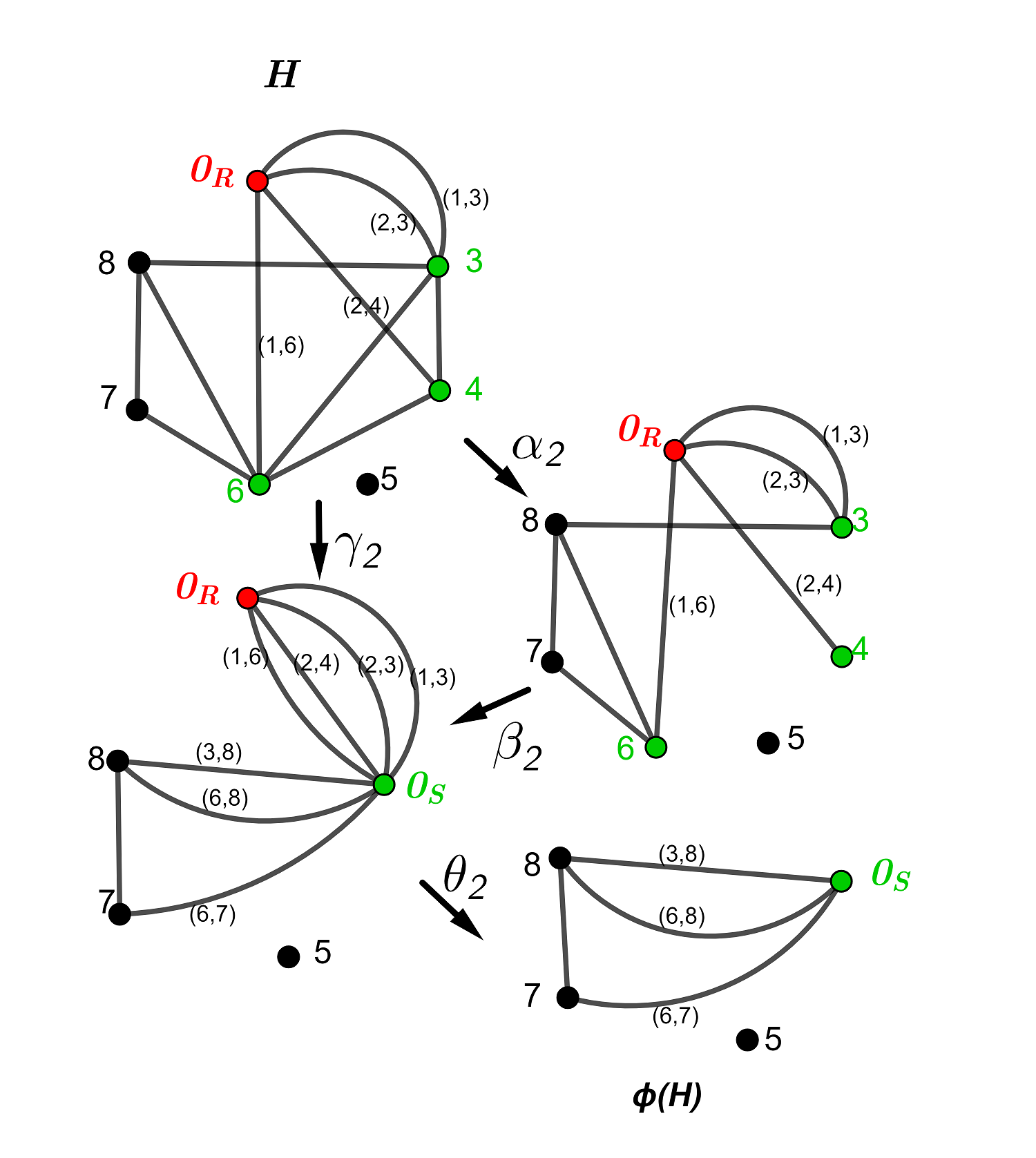}
\end{center}
\caption{The successive images of a graph $H\in \mathcal{H}(K_{n/r},S)$ for $%
\protect\phi =\protect\theta _{2}\circ \protect\beta _{2}\circ \protect\alpha
_{2}$, $R=\{1,2\}$ and $S=\{3,4,6\}$.}
\end{figure}

Now, starting from an arbitray spanning subgraph $H^{\prime }$ of $K_{n-r/s}$%
, the inverse image $\theta _{2}^{-1}\left( H^{\prime }\right) $ is obtained
by adding:

$l_{1}$ edges among the edges labeled $\left\{ u_{i},v_{1}\right\} $, with $%
u_{i}\in R$ and $1\leq l_{1}\leq r$,

$l_{2}$ edges among the edges labeled $\left\{ u_{i},v_{2}\right\} $, with $%
u_{i}\in R$ and $1\leq l_{2}\leq r$,

...

$l_{s}$ edges among the edges labeled $\left\{ u_{i},v_{s}\right\} $, with $%
u_{i}\in R$ and $1\leq l_{s}\leq r$.

Once $l_{1},l_{2},...,l_{S}$ are fixed, we have, for each of the graphs $%
H^{0}\in \theta _{2}^{-1}\left( H^{\prime }\right) $ and $H^{\prime }$, the
relations: 
\begin{equation*}
v\left( H^{0}\right) =v\left( H^{\prime }\right) +1,\;e\left( H^{0}\right)
=e\left( H^{\prime }\right) +l_{1}+l_{2}+...+l_{s},\;\;c\left( H^{0}\right)
=c\left( H^{\prime }\right) .
\end{equation*}
Using equation $\left( 3.5\right) $, we obtain the following relations
between $H^{\ast }=\beta _{2}^{-1}\left( H^{0}\right) $ and $H^{\prime }$: 
\begin{equation*}
v\left( H^{\ast }\right) =v\left( H^{\prime }\right) +s,\;e\left( H^{\ast
}\right) =e\left( H^{\prime }\right) +l_{1}+l_{2}+...+l_{s},\;\;c\left(
H^{\ast }\right) =c\left( H^{\prime }\right) .
\end{equation*}
By summing over all graphs $H^{\ast }\in \alpha _{2}\left( \mathcal{H}\left(
K_{n/r},S\right) \right) $, it follows that: 
\begin{equation*}
\sum\limits_{H^{\ast }}\left( x-1\right) \left( y-1\right) ^{e\left( H^{\ast
}\right) +c\left( H^{\ast }\right) -v\left( H^{\ast }\right)
}=\prod\limits_{j=1}^{s}\sum\limits_{l_{j}=1}^{r}\binom{r}{l_{j}}\left(
y-1\right) ^{l_{j}-1}\sum\limits_{H^{\prime }}\left( x-1\right) ^{c\left(
H^{\prime }\right) -1}\left( y-1\right) ^{e\left( H^{\prime }\right)
+c\left( H^{\prime }\right) -v\left( H^{\prime }\right) }\text{,}
\end{equation*}
where the last sum runs over all spanning subgraph $H^{\prime }$\ of $%
K_{n-r/s}.$ We therefore obtain, using Eq. $\left( 3.4\right) $: 
\begin{equation*}
\sum\limits_{H^{\ast }}\left( x-1\right) \left( y-1\right) ^{e\left( H^{\ast
}\right) +c\left( H^{\ast }\right) -v\left( H^{\ast }\right) }=\left( \left[
r\right] _{y}\right) ^{s}T_{n-r}^{\left( s\right) }(x,y)\text{.}
\end{equation*}
It remains to enumerate the graphs $H\in \mathcal{H}\left( K_{n/r},S\right) $
such that $H\in \alpha _{2}^{-1}\left( H^{\ast }\right) $. To do this, one
can add to $H^{\ast }$, $p$ edges chosen among the edges of $K_{S}$ with $%
0\leq p\leq \binom{s}{2}$. If $p$ is fixed, there are $\binom{\binom{s}{2}}{p%
}$ possible choices for these $p$ edges. For each choice, we have the
relations: 
\begin{equation*}
v\left( H\right) =v\left( H^{\ast }\right) ,\;e\left( H\right) =e\left(
H^{\ast }\right) +p,\;c\left( H\right) =c\left( H^{\ast }\right) \text{,}
\end{equation*}
where the last equality holds because all elements of $S$ belong to the same
connected component containing $0_{R}$. It follows that: 
\begin{eqnarray*}
T_{S}\left( x,y\right)  &=&\sum\limits_{p=0}^{\binom{s}{2}}\binom{\binom{s}{2%
}}{p}\left( y-1\right) ^{p}\sum\limits_{H^{\ast }}\left( x-1\right) \left(
y-1\right) ^{e\left( H^{\ast }\right) +c\left( H^{\ast }\right) -v\left(
H^{\ast }\right) } \\
&=&y^{\binom{s}{2}}\left( \left[ r\right] _{y}\right) ^{s}T_{n-r}^{\left(
s\right) }\left( x,y\right) \text{.}
\end{eqnarray*}
By adding all possible choice of $S$, with $s=$ $\left| S\right| $ between 2
and $n-r$, the contribution of case 3) is therefore: 
\begin{equation*}
\sum\limits_{s=2}^{n-r}\binom{n-r}{s}y^{\binom{s}{2}}\left[ r\right]
_{y}^{s}T_{n-r}^{\left( s\right) }\left( x,y\right) \text{.}
\end{equation*}
By adding the contributions of 1), 2) and 3), we indeed obtain Equation $%
\left( 1.3\right) $.
\end{proof}

\begin{remark}
It is easy to see that one could also write $\phi =\beta _{1}\circ \alpha
_{1}\circ \theta _{1}$ by defining $\alpha _{1}$ and $\beta _{1}$ in a
manner analogous to $\alpha _{2}$ and $\beta _{2}$. However, this
decomposition does not allow for an easy computation of the contribution, in
the sum of \ Equation $\left( 3.3\right) $, corresponding to the preimage $%
\phi ^{-1}\left( H^{\prime }\right) $ of a spanning subgraph of $%
K_{V-R}/E\left( K_{S}\right) \approx K_{n-r/s}$. Indeed, unlike $\alpha _{2}$
and $\beta _{2}$, the maps $\alpha _{1}$ and $\beta _{1}$ do not preserve
the number of connected components. This is illustrated in Figure 3, which
shows the successive images in the decomposition $\phi =\beta _{1}\circ
\alpha _{1}\circ \theta _{1}$ of the same graph $H$ as in Figure 2. In fact,
the variation in the number of connected components under $\alpha _{1}$ and $%
\beta _{1}$ depends of the graph $H$, which makes the computation
impraticable with this decomposition.
\end{remark}

\bigskip 
\begin{figure}[!htp]
\begin{center}
\includegraphics[width=12cm]{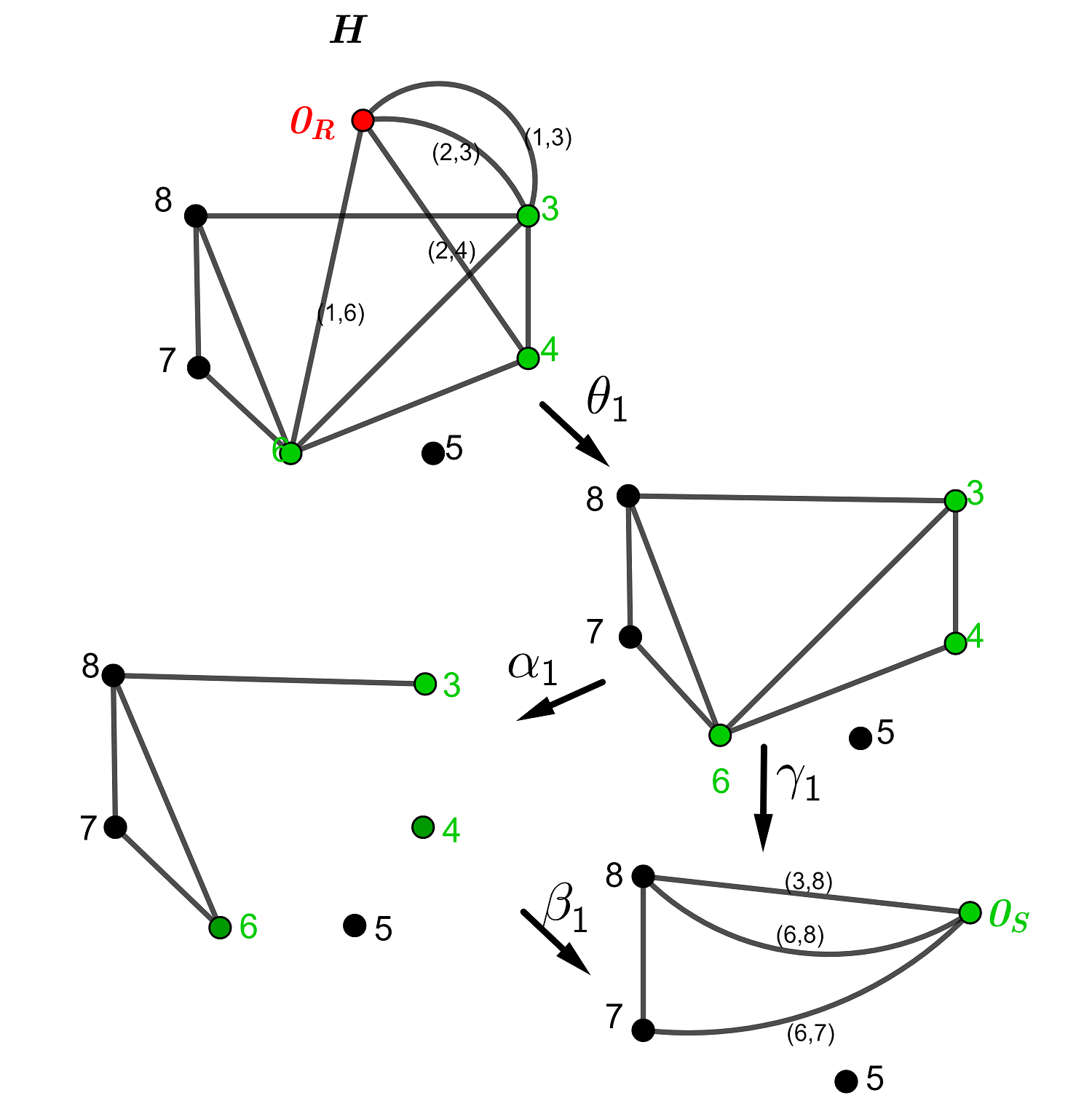}
\end{center}
\caption{The successive images of the graph $H$ from Figure 2 for $\protect%
\phi=\protect\beta_1\circ\protect\alpha_1\circ\protect\theta_1$, $R=\{1,2\}$
and $S=\{3,4,6\}$.}
\end{figure}

The graph $K_{n/n}$ is simply the graph with a single vertex, $0_{R}$, and
no edges,\ whose Tutte polynomial is $T_{n}^{\left( n\right) }\left(
x,y\right) =1$. As we did for the $J_{n}^{\left( r\right) }$ in [2], let us
place the polynomials $T_{n}^{\left( r\right) }$ in a triangular array whose
rows and columns are indexed by $n$ and $r$, respectively. Equation $\left(
1.3\right) $ allows us to recursively calculate row by row all the
polynomials $T_{n}^{\left( r\right) }$ for $n>r\geq 1$. In particular, the
Tutte polynomials of the complete graphs $K_{n}$, constitute the first
column of this array. This triangular array is represented below for $1\leq
r\leq n\leq 5$.

\bigskip

\begin{tabular}{|l|l|l|l|l|l|}
\hline
\begin{tabular}{ll}
& $r$ \\ 
$n$ & 
\end{tabular}
& $1$ & $2$ & $3$ & $4$ & $5$ \\ \hline
$1$ & $1$ &  &  &  &  \\ \hline
$2$ & $x$ & $1$ &  &  &  \\ \hline
$3$ & $x+y+x^{2}$ & $x+y$ & $1$ &  &  \\ \hline
$4$ & 
\begin{tabular}{l}
$2x+2y+3x^{2}+4xy$ \\ 
$+3y^{2}+x^{3}+y^{3}$%
\end{tabular}
& 
\begin{tabular}{l}
$x+y+x^{2}+2xy$ \\ 
$+2y^{2}+y^{3}$%
\end{tabular}
& $x+y+y^{2}$ & $1$ &  \\ \hline
$5$ & $
\begin{tabular}{l}
$6x+6y+11x^{2}+20xy$ \\ 
$+15y^{2}+6x^{3}+10x^{2}y$ \\ 
$15xy^{2}+15y^{3}+x^{4}+5xy^{3}$ \\ 
$+10y^{4}+5y^{5}+y^{6}$%
\end{tabular}
$ & $
\begin{tabular}{l}
$2x+2y+3x^{2}+7xy$ \\ 
$+6y^{2}+x^{3}+3x^{2}y$ \\ 
$+6xy^{2}+7y^{3}+3xy^{3}$ \\ 
$+6y^{4}+3y^{5}+y^{6}$%
\end{tabular}
$ & 
\begin{tabular}{l}
$x+y+x^{2}+2xy$ \\ 
$+2y^{2}+2xy^{2}+3y^{3}$ \\ 
$+2y^{4}+y^{5}$%
\end{tabular}
& $
\begin{tabular}{l}
$x+y$ \\ 
$+y^{2}+y^{3}$%
\end{tabular}
$ & $1$ \\ \hline
\end{tabular}

\medskip \textbf{Table 1: Polynomials }$T_{n}^{\left( r\right) }$\textbf{\
for }$5\geq n\geq r\geq 1$.\newline

\section{\protect\bigskip Case of connected graphs}

We have seen that Eq. $\left( 1.1\right) $ is the special case of Eq. $%
\left( 1.3\right) $ obtained for $x=1$ and $y=q$. Moreover, for any
connected graph $G$ such that $\left| V\left( G\right) \right| =n$, we
deduce from the definition of the Tutte polynomial used in the proof of
Theorem 1.1 that: 
\begin{equation}
T_{G}\left( 1,y\right) =\sum\limits_{C}\left( y-1\right) ^{e\left( C\right)
+1-n}\text{,}  \tag{4.1}
\end{equation}
where the sum is over the spanning connected subgraphs of $G$. By applying
this formula, Eq. $\left( 1.1\right) $ could be directly proved in a way
similar to the proof of Theorem 1.1, using the decomposition: 
\begin{equation}
\mathcal{C}\left( K_{n/r}\right) =\biguplus\limits_{S\subseteq \overline{R}}%
\mathcal{C}\left( K_{n/r},S\right)
=\biguplus\limits_{s=0}^{n-r}\biguplus_{\left| S\right| =s}\mathcal{C}\left(
K_{n/r},S\right) \text{.}  \tag{4.2}
\end{equation}
In this equation, $\mathcal{C}\left( K_{n/r}\right) $ is the set of spanning
connected subgaphs of $K_{n/r}$ and $\mathcal{C}\left( K_{n/r},S\right) $ is
the subset of $\mathcal{C}\left( K_{n/r}\right) $ whose graphs have exactly
the subset $S$ as vertices adjacent to $0_{R}$. 

Let us define: 
\begin{equation}
C_{n}^{\left( r\right) }(t)=\sum\limits_{C\in \mathcal{C}\left(
K_{n/r}\right) }t^{e\left( C\right) }\text{,}  \tag{4.3}
\end{equation}
the relation $\left( 4.1\right) $ applied to $K_{n/r}$, shows that, taking $%
\left( 1.2\right) $ into account: 
\begin{equation*}
J_{n}^{\left( r\right) }\left( q\right) =\left( q-1\right) ^{-\left(
n-r\right) }\sum\limits_{C\in \mathcal{C}\left( K_{n/r}\right) }\left(
q-1\right) ^{e\left( C\right) }\text{.}
\end{equation*}
By comparing with $\left( 4.3\right) $, we obtain: 
\begin{equation}
C_{n}^{\left( r\right) }(t)=t^{n-r}J_{n}^{\left( r\right) }\left( 1+t\right) 
\text{.}  \tag{4.4}
\end{equation}
Equation $\left( 4.4\right) $ generalizes the following case for $r=1$: 
\begin{equation}
C_{n}\left( t\right) =t^{n-1}J_{n}\left( 1+t\right) \text{,}  \tag{4.5}
\end{equation}
which was first proven algebraically in [6] and later combinatorially in
[5]. Note also that $\left( 4.4\right) $ can be obtained as a special case
of Corollary $7.1$ of [8], up to changes in notation.

\medskip

\textbf{References}

$\left[ 1\right] $J.A. Bondy, U.S.R Murty, \textit{Graph Theory}. Springer,
London (2008).

$\left[ 2\right] $V. Brugidou, A q-analog of certain symmetric functions and
one of its specializations, (2023), arXiv:2302.11221.

$\left[ 3\right] $ V. Brugidou, On a particular specialization of monomial
symmetric functions, (2023), preprint, arXiv:2306.15300.

$\left[ 4\right] $V. Brugidou, On a particular specialization of monomial
symmetric functions, International Conference on Enumerative Combinatorics
and Applications (ICECA), (September 4-6-, 2023).

$\left[ 5\right] $ I. M. Gessel and D. L. Wang, Depth-first search as a
combinatorial correspondence. \textit{J. Combin. Theory Ser. A}, \textbf{26}
(1979), 308-313.

$\left[ 6\right] $ C. L. Mallows and J. Riordan, The inversion enumerator
for labeled trees. \textit{Bull. Amer. Soc. }\textbf{74} (1968) 92-94.

$\left[ 7\right] $ R.P. Stanley, Hyperplane arrangements, parking functions,
and tree inversions. In\textit{\ Mathematical Essays in Honor of Gian-carlo
Rota} \textit{(Cambridge, MA, 1996)}, vol. 161 of Progr. Math.,
Birkh\"{a}user, Boston, MA (1998) 359-375.

$\left[ 8\right] $ C. H. Yan. Generalized parking functions, tree
inversions, and multicolored graphs. \textit{Adv. in Appl. Math}. \textbf{27}%
(2-3) (2001) 641-670.

$\left[ 9\right] $ C.H. Yan, Parking functions, in M. Bona (ed.), \textit{%
Handbook of Enumerative Combinatorics}, CRC Press, Boca Raton, FL (2015) pp.
835-893.

\bigskip

\end{document}